# Eigenvalues of an axially loaded cantilever beam with an eccentric end rigid body

**Seyed Amir Mousavi Lajimi**
Faculty of Engineering, University of Waterloo
200 University Avenue West
Waterloo, Ontario, Canada, N2L 3G1
`samousavilajimi@uwaterloo.ca, s.a.m.lajimi@a3.epfl.ch`

**G. R. Heppler**
Faculty of Engineering, University of Waterloo
200 University Avenue West
Waterloo, Ontario, Canada, N2L 3G1
`heppler@uwaterloo.ca`

***Abstract*** - An analytical form of the characteristic equation for a vertically mounted cantilever beam with an end rigid body is obtained and solved for the eigenvalues of the structure. The effect of the weight of the structure is taken into consideration by estimating the load as a function of the length of the beam. The mass, rotary inertia and eccentricity of the end rigid body are demonstrated to considerably affect the eigenvalues of the structure.

## 1. Introduction

This work concerns the eigenvalues of vertical cantilever beams under self-weight carrying an eccentric rigid end body. Self-weight can be observed as a linearly varying axial load representing a non-follower force. Because of the apparent difficulty in solving the corresponding differential eigenvalue problem of a vertical cantilever beam with an eccentric end rigid body, ordinarily the problem has been simplified in some senses in the past to obtain solutions. Therefore, researchers have mainly used approximation methods to estimate eigenvalues and eigenfunctions of even simpler cases.

A simplified case, where the varying load was approximated by a constant load, was considered by Bokaian (1998) and the frequencies of a uniform beam under constant axial compressive loads were computed. Later, Bokaian extended his work to various combinations of conventional boundary conditions (Bokaian, 1990). A list of earlier references is given in Bokaian's work which gives an understanding of the long history of the subject. Bokaian's works did not include any boundary conditions with an end mass. In contrast, Naguleswaran (2006) considered nonclassical boundary conditions including end rigid bodies and inertias of end bodies, but no axial force was included.

Earlier, Paidoussis and Des Trois Maissons studied the free flexural vibration of an internally damped cantilever including the self-weight of the beam (Paidoussis and Des Trois Maissons, 1971). They used the Galerkin-type approximation method where the basis functions were the eigenfunctions of a horizontal cantilevered beam. Naguleswaran studied the effects of a linearly varying axial force on the natural frequencies of a uniform single span beam with classical or conventional boundary conditions by using Frobenius power series solutions (Naguleswaran 1991 and 2004). Using the Rayleigh-Ritz method, Schafer (Schafer 1985) studied the effect of gravity on the natural frequencies and mode shapes of a hanging cantilever beam.

In this work, we consider a general case where the vertically standing beam is influenced by its self-weight and a general nonclassical boundary condition, *i.e.* an end rigid body with an eccentric center of mass. The mass moment of inertia of the end rigid body is also taken into consideration. The differential eigenvalue problem is solved by using the Frobenius method and the eigenvalues of the structure are obtained and presented in the form of tables and plots.

## 2. Problem Formulation

A schematic of the uniform cantilever beam carrying an eccentric end mass is presented in Fig. 1. The positive Z-axis of the coordinate system points upward while the origin is attached to the ground.

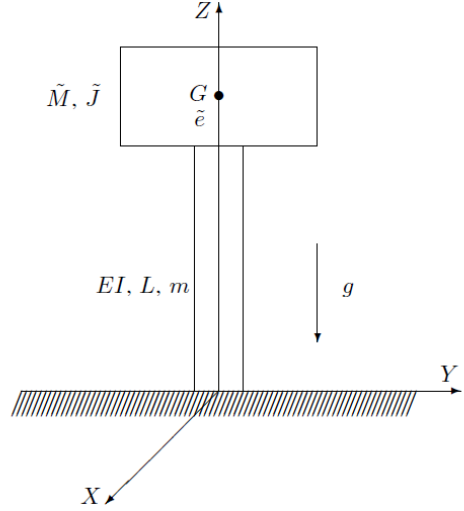

Figure 1. A schematic diagram of the clamped standing beam carrying an end rigid body

### 2. 1. Differential Equation of Motion

The beam has a constant cross-section area $A$, moment of inertia $I$, length $L$, mass per unit length $m$, and modulus of elasticity $E$. The cross-sectional area of the beam is uniform and its material is homogenous. The parameters $\tilde{M}$, $G$, $\tilde{J}$, and $\tilde{e}$ represent the mass, the center of mass, the rotary inertia, and the center of mass off-set values of the end mass. The principles of Newtonian mechanics are employed to derive the partial differential equation of motion and the boundary conditions. Considering a small element of the beam, the equilibrium equation in a transverse (horizontal) direction is written as

$$\frac{\partial S}{\partial Z} - C\frac{\partial V}{\partial \tau} + Q = m\frac{\partial^2 V}{\partial \tau^2} \qquad (1)$$

where $S$ indicates the shear force, $Q$ the transverse force, $C$ the damping factor, and $V$ the transverse displacement component of the beam. The moment balance equation for the beam element becomes

$$S - \frac{\partial M}{\partial Z} + P(Z,\tau)\frac{\partial V}{\partial Z} = 0 \qquad (2)$$

on account of a small rotation assumption which means $\sin\frac{\partial V}{\partial Z} \approx \frac{\partial V}{\partial Z}$, and $P(Z,\tau)$ represents the axially varying force. Substituting Eq. (2) into Eq. (1) and using the appropriate relation for the moment in terms of elastic rigidity, gives the differential equation of motion in its final form as

$$\frac{\partial^2}{\partial Z^2}\left(EI\frac{\partial^2 V}{\partial Z^2}\right) + \frac{\partial}{\partial Z}\left(P(Z,\tau)\frac{\partial V}{\partial Z}\right) + C\frac{\partial V}{\partial \tau} + m\frac{\partial^2 V}{\partial \tau^2} = Q \qquad (3)$$

The boundary conditions of the beam at $Z = 0$ are readily seen to be

$$V = 0 \quad \text{and} \quad \frac{\partial V}{\partial Z} = 0 \qquad (4)$$

To find the boundary condition at the other end of the beam the acceleration of the rigid end body should be computed. The transverse acceleration of the centre of mass of the end rigid body is

$$a_G = -\frac{\partial^2 V}{\partial \tau^2} - \tilde{e}\frac{\partial^2}{\partial \tau^2}\frac{\partial V}{\partial Z} \tag{5}$$

A shear force balance in the transverse direction gives:

$$S = \tilde{M}\frac{\partial^2 V}{\partial \tau^2} + \tilde{e}\,\tilde{M}\frac{\partial^2}{\partial \tau^2}\frac{\partial V}{\partial Z} \tag{6}$$

Then, the boundary condition at $Z = L$ is obtained as:

$$\frac{\partial}{\partial Z}\left(EI\frac{\partial^2 V}{\partial Z^2}\right) + P(Z,\tau)\frac{\partial V}{\partial Z} - \tilde{M}\frac{\partial^2 V}{\partial \tau^2} - \tilde{e}\,\tilde{M}\frac{\partial^2}{\partial \tau^2}\frac{\partial V}{\partial Z} = 0 \tag{7}$$

The second natural boundary condition is obtained using the definition of angular momentum of the end rigid body and moment balance about the $Z = L$ end of the beam as

$$\tilde{M}\,\tilde{e}\,\frac{\partial^2 V}{\partial \tau^2} + (\tilde{J} + \tilde{M}\,\tilde{e}^2)\frac{\partial^2}{\partial \tau^2}\frac{\partial V}{\partial Z} + EI\frac{\partial^2 V}{\partial Z^2} = 0 \tag{8}$$

Therefore, the governing differential equation and boundary conditions are given in Eqs. (3), (4), (7) and (8). By representing the model in dimensionless form, the number of parameters in the equation of motion and boundary conditions is reduced. To transform the system from the dimensional to the dimensionless form, the following set of variables is introduced to the governing differential equations and boundary conditions:

$$z = \frac{Z}{L}, \qquad v = \frac{V}{\ell}, \qquad t = \frac{\tau}{\varkappa} \tag{9}$$

where $\ell$ and $\varkappa$ represent two arbitrary constants in time and length units. Consequently, the equation of motion and boundary conditions for an axially loaded cantilever beam with an end rigid mass are

$$\frac{\partial^4 v}{\partial z^4} + \frac{\partial}{\partial z}\left(p(z)\frac{\partial v}{\partial z}\right) + c\frac{\partial v}{\partial t} + \frac{\partial^2 v}{\partial t^2} = q \tag{10}$$

$$v = 0, \qquad at\ z = 0 \tag{11}$$

$$\frac{\partial v}{\partial z} = 0, \qquad at\ z = 0 \tag{12}$$

$$\frac{\partial^3 v}{\partial z^3} + p(z)\frac{\partial v}{\partial z} - M\frac{\partial^2 v}{\partial t^2} - M\,e\,\frac{\partial^2}{\partial t^2}\frac{\partial v}{\partial z} = 0, \qquad at\ z = 1 \tag{13}$$

$$M\,e\,\frac{\partial^2 v}{\partial t^2} + (J + M\,e^2)\frac{\partial^2}{\partial t^2}\frac{\partial v}{\partial z} + \frac{\partial^2 v}{\partial z^2} = 0, \qquad at\ z = 1 \tag{14}$$

where

$$p(z) = \frac{P(Z)\,L^2}{EI},\ c = \frac{C\,L^4}{EI\,\varkappa},\ q = \frac{Q\,L^4}{EI\,\ell},\ e = \frac{\tilde{e}}{L},\ M = \frac{\tilde{M}\,L^3}{EI\,\varkappa},$$

$$J = \frac{\tilde{J}\,L}{EI\,\varkappa^3},\qquad \ell = L,\qquad \varkappa = \sqrt{\frac{EI}{m\,L^2}} \tag{15}$$

are the parameters of the system.

## 2. 2. Eigenvalue Problem and the Method of Solution

To define the eigenvalue problem the external and damping forces are set to zero, $q = 0$ and $c = 0$, and a harmonic solution in time, $v(z,t) = \eta(z)\,\sin(\sqrt{\Lambda}\,t)$, is substituted into Eqs. (10) – (14) to obtain the following differential equation and boundary conditions:

$$\frac{d^4\eta}{dz^4} + \frac{d}{dz}\left(p(z)\frac{d\eta}{dz}\right) = \Lambda\,\eta \tag{16}$$

$$\eta = 0 \text{ and } \frac{d\eta}{dz} = 0, \quad at\ z = 0 \tag{17}$$

$$\frac{d^3\eta}{dz^3} + p(z)\frac{d\eta}{dz} + \Lambda\,M\left(\eta + e\frac{d\eta}{dz}\right) = 0, \quad at\ z = 1 \tag{18}$$

$$\frac{d^2\eta}{dz^2} - \Lambda\left(M\,e\,\eta + (J + M\,e^2)\frac{d\eta}{dz}\right) = 0, \quad at\ z = 1 \tag{19}$$

which govern the pattern of motion for linear modes of the clamped beam carrying an eccentric end rigid body. Because the system features an eigenvalue problem with variable coefficients in the variable $z$, i.e. $p(z)$, the Frobenius method is employed to solve for the eigenvalues and eigenfunctions of the system (Hildebrand 1976). To this end, the spatial function, $\eta(z)$, is expressed as a power series such that

$$\eta_r(z) = \sum_{n=0}^{\infty} a_n(r) z^{n+r} \tag{20}$$

where $r$ is an undetermined variable which is obtained from the corresponding indicial equation of the differential equation. Substituting Eq.(20) into Eq.(16), and noting that for the problem under investigation

$$p(z) = p_0 + \gamma\,z \tag{21}$$

where the slope of the axial load function, $p(z)$, given by

$$\gamma = (p_1 - p_0)/(z_1 - z_0) \tag{22}$$

renders the differential eigenvalue problem into the power series form

$$\begin{aligned}
&r(r-1)(r-2)(r-3)a_0(r)z^{r-4} + (r+1)r(r-1)(r-2)a_1(r)z^{r-3} \\
&+ (r+2)(r+1)r(r-1)a_2(r)z^{r-2} \\
&+ (r+3)(r+2)(r+1)r a_3(r)z^{r-1} \\
&+ p_0\,r(r-1)a_0(r)z^{r-2} + p_0\,(r+1)r\,a_1(r)z^{r-1} \\
&+ \gamma\,r^2 a_0(r)z^{r-1} \\
&+ \sum_{n=4}^{\infty}[(n+r)(n+r-1)(n+r-2)(n+r-3)a_n(r) \\
&+ p_0\,(n+r-3)(n+r-2)a_{n-2}(r) \\
&+ \gamma(n+r-3)^2 a_{n-3}(r) - \Lambda a_{n-4}(r)]\,z^{n+r-4} = 0
\end{aligned} \tag{23}$$

The coefficients of the $z^{n+r}$ in Eq.(23) must simultaneously go to zero producing four series functions of the length variable $z$ corresponding to each root of the indicial equation, i.e. $r = 0, 1, 2,$ and $3$. The general solution of the differential equation is then the sum of all solutions:

$$\eta(z) = \sum_{r=0}^{3} A(r)\eta_r(z) \tag{24}$$

where the unknown coefficients, $A(r)$, are obtained by applying the boundary conditions. The frequency equation is then obtained by setting the determinant of the coefficients matrix to zero resulting in

$$\begin{aligned}
&\{\eta_3'''(z) - (\gamma + p_0)\,\eta_3'(z) + \Lambda\,M\,[\eta_3(z) + e\,\eta_3'(z)]\}\,\{\eta_2''(z) \\
&\quad - \Lambda\,[M\,e\,\eta_2(z) + (J + M\,e^2)\,\eta_2'(z)]\} \\
&- \{\eta_2'''(z) - (\gamma + p_0)\,\eta_2'(z) \\
&\quad + \Lambda\,M\,[\eta_2(z) + e\,\eta_2'(z)]\}\,\{\eta_3''(z) \\
&\quad - \Lambda\,[M\,e\,\eta_3(z) + (J + M\,e^2)\,\eta_3'(z)]\} = 0
\end{aligned} \tag{25}$$

where the $\eta_r(z)$ terms indicate the power series including $\Lambda$ as the coefficient in terms higher than fourth order. The frequency equation, Eq.(25), is numerically solved for the natural frequencies of the structure.

## 3. Results and Discussions

To find the roots of the characteristic equation, Eq.(25), the sign of the left hand side of the characteristic equation is evaluated for a range of eigenvalues starting from zero with an assumed increment size in the eigenvalue. The process is repeated with different increment sizes to avoid missing roots. Sign alterations indicate intervals where roots of the characteristic equation should be searched for. Then, implementing a bisection algorithm the root of the characteristic equation in each interval is found with a preset accuracy.

Table 1. Basic structural specifications based on (Cai and Chen 1996)

| Description | Numerical value |
|---|---|
| Modulus of elasticity | 7.9336 GPa |
| Linear mass density | 61.08 kgm$^{-1}$ |
| Outer diameter | 1.238 m |
| Inner diameter | 1.219 m |
| Gravity acceleration | 9.81 ms$^{-2}$ |

To determine the number of terms required in the series, the axial load is set to zero once and the results of this study are compared with the exact solution computed in (Mousavi Lajimi and Heppler 2012). The input data is given in Table 1 based on the parameter values used by (Cai and Chen 1996). Figures 2 and 3 show the variation in the first and second eigenvalue (natural frequency squared) as a function of the length of the structure which is proportional to the weight of the beam. To identify the critical length, which corresponds to the self-weight buckling of the structure, the zero-crossing of the eigenvalue curve in Fig. 2 should be obtained.

In Figs. 4 and 5 the evolution of the fundamental eigenvalues as a function of the mass of the end rigid body are presented. These two plots correspond to the cantilever beam with a point mass at the end of the beam, *i.e.* $J = 0$ and $e = 0$. In Tables 2 and 3 the effects of increasing mass moment of inertia of the end rigid body to the eigenvalues are presented. Table 2 corresponds the case of a thin end rigid body, *i.e.* $e = 0$, and Table 3 represents the case of nonzero eccentricity. Finally, Figs. 6 and 7 demonstrate how increasing the eccentricity influences the first two eigenvalues of the system.

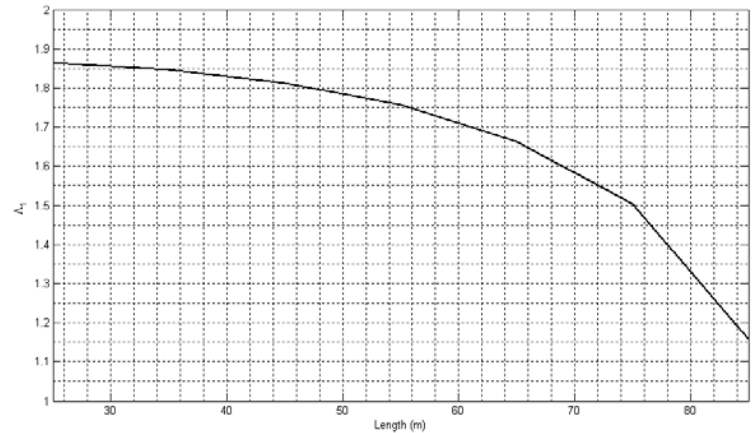

Figure 2. The lowest eigenvalue versus length of the beam with no tip mass, *i.e.* $M = 0, J = 0$ and $e = 0$

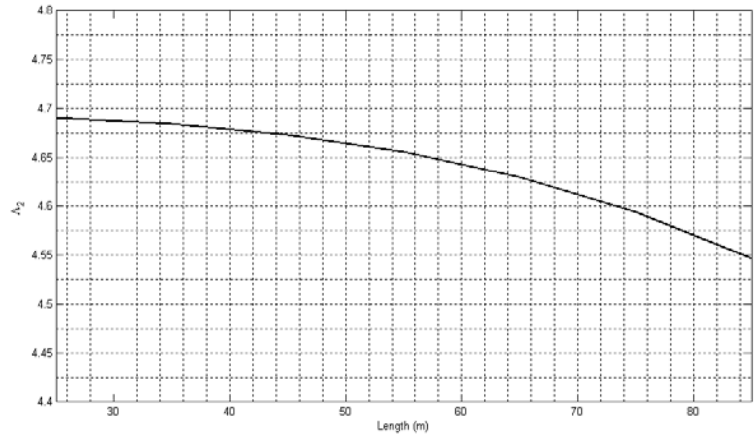

Figure 3. The second lowest eigenvalue versus length of the beam with no tip mass, *i.e.* $M = 0, J = 0$ and $e = 0$

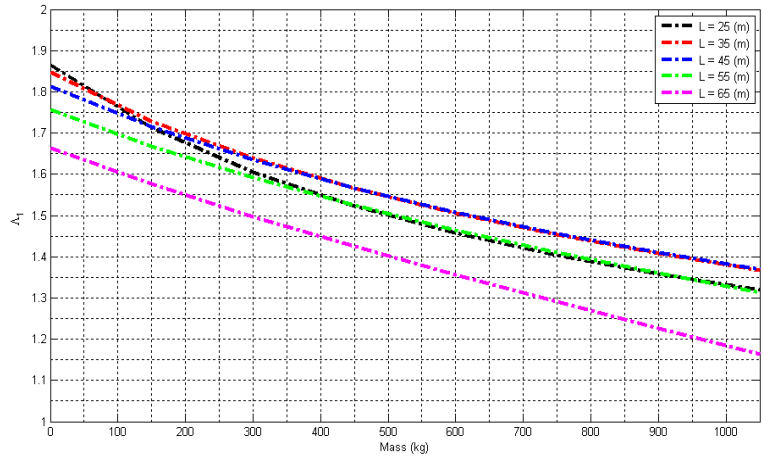

Figure 4. The lowest eigenvalue versus mass of end rigid body for $J = 0$ and $e = 0$

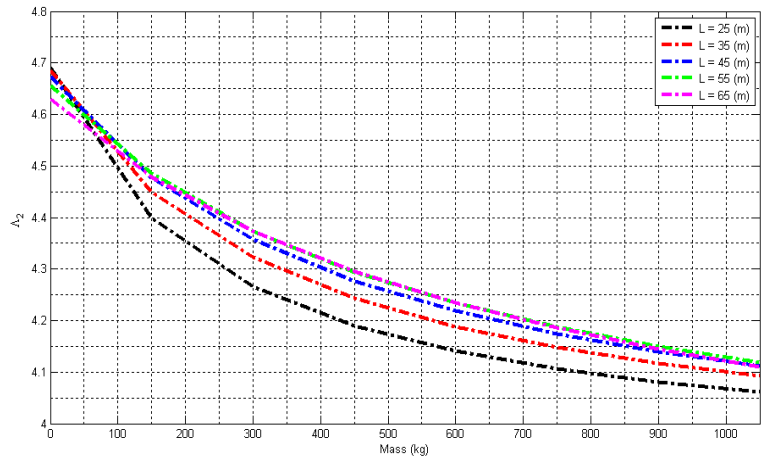

Figure 5. The lowest eigenvalue versus mass of end rigid body for $J = 0$ and $e = 0$

Table 2. Variation in the first and second eigenvalue for $\widetilde{M} = 600$kg and $\tilde{e} = 0$m

| $L$ (m) | Eigenvalue | $\tilde{J}$ | 0 | 5000 | 10000 | 15000 | 20000 | 25000 |
|---|---|---|---|---|---|---|---|---|
| 25 | $\Lambda_1$ | | 1.458 | 1.452 | 1.446 | 1.439 | 1.433 | 1.427 |
|  | $\Lambda_2$ | | 4.140 | 3.875 | 3.646 | 3.464 | 3.32 | 3.202 |
| 35 | $\Lambda_1$ | | 1.505 | 1.502 | 1.499 | 1.497 | 1.494 | 1.491 |
|  | $\Lambda_2$ | | 4.188 | 4.079 | 3.972 | 3.871 | 3.778 | 3.692 |
| 45 | $\Lambda_1$ | | 1.508 | 1.506 | 1.505 | 1.503 | 1.502 | 1.500 |
|  | $\Lambda_2$ | | 4.219 | 4.164 | 4.108 | 4.054 | 4.000 | 3.949 |
| 55 | $\Lambda_1$ | | 1.465 | 1.464 | 1.463 | 1.463 | 1.462 | 1.461 |
|  | $\Lambda_2$ | | 4.234 | 4.202 | 4.17 | 4.139 | 4.107 | 4.076 |
| 65 | $\Lambda_1$ | | 1.357 | 1.356 | 1.356 | 1.355 | 1.355 | 1.354 |
|  | $\Lambda_2$ | | 4.234 | 4.214 | 4.194 | 4.174 | 4.154 | 4.134 |

Table 3. Variation in the first and second eigenvalue for $\widetilde{M} = 600$kg and $\tilde{e} = 8$m

| $L$ (m) | Eigenvalue | $\tilde{J}$ | 0 | 5000 | 10000 | 15000 | 20000 | 25000 |
|---|---|---|---|---|---|---|---|---|
| 25 | $\Lambda_1$ | | 1.266 | 1.262 | 1.258 | 1.255 | 1.251 | 1.247 |
|  | $\Lambda_2$ | | 3.280 | 3.215 | 3.157 | 3.104 | 3.055 | 3.011 |
| 35 | $\Lambda_1$ | | 1.376 | 1.373 | 1.371 | 1.369 | 1.367 | 1.365 |
|  | $\Lambda_2$ | | 3.509 | 3.466 | 3.426 | 3.388 | 3.352 | 3.317 |
| 45 | $\Lambda_1$ | | 1.417 | 1.415 | 1.414 | 1.413 | 1.412 | 1.411 |
|  | $\Lambda_2$ | | 3.678 | 3.649 | 3.621 | 3.594 | 3.568 | 3.542 |
| 55 | $\Lambda_1$ | | 1.399 | 1.399 | 1.398 | 1.397 | 1.397 | 1.396 |
|  | $\Lambda_2$ | | 3.800 | 3.780 | 3.760 | 3.741 | 3.722 | 3.703 |
| 65 | $\Lambda_1$ | | 1.310 | 1.310 | 1.309 | 1.309 | 1.308 | 1.308 |
|  | $\Lambda_2$ | | 3.881 | 3.867 | 3.853 | 3.839 | 3.825 | 3.812 |

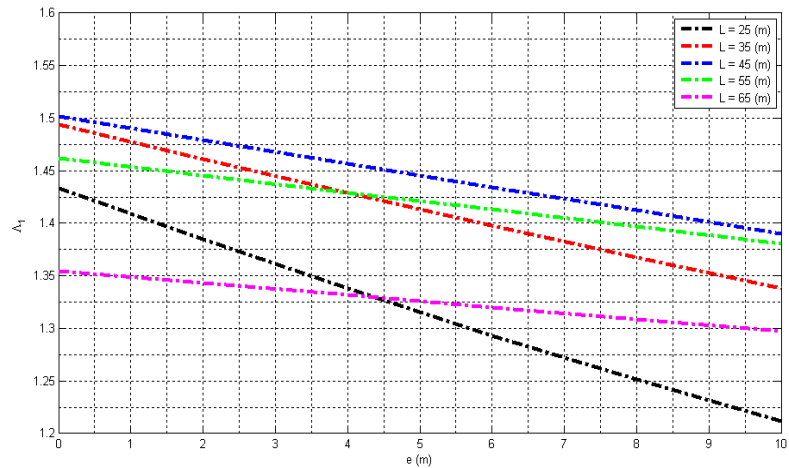

Figure 6. The lowest eigenvalue versus eccentricity of the center of mass of end rigid body for $\tilde{J} = 20000$ kg m$^2$ and $\widetilde{M} = 600$ kg

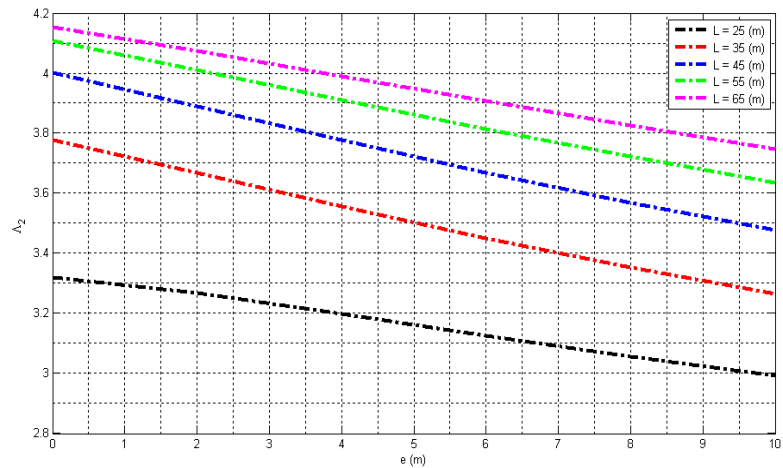

Figure 7. The lowest eigenvalue versus eccentricity of the center of mass of end rigid body for $\tilde{J} = 20000$ kg m$^2$ and $\widetilde{M} = 600$ kg

## 4. Conclusion

The eigenvalues of a cantilever beam with an eccentric end rigid body have been studied by using the method of Frobenius for solving differential equations with variable coefficients. It has been shown that all parameters appear in the characteristic equation and therefore influence the eigenvalues of the structure. The length of the beam representing the self-weight of the structure as well as the mass of the end rigid body have the largest effect in reducing the eigenvalues of the structure. As the moment of inertia of the end rigid body becomes larger, the eigenvalues of the structure are reduced; however, in compare with the variation in the eccentricity, the magnitude of the moment of inertia is less effective in changing the eigenvalues of the structure.